\documentclass[11pt]{article}
\usepackage{amsmath}
\usepackage{amsmath,amssymb,amsfonts,amsthm,fancyhdr}
\usepackage{epsfig,graphicx,picins,picinpar,subfigure}
\usepackage{pstricks}
\usepackage{fancyvrb}
\usepackage[numbers,sort&compress]{natbib}

\begin{document}

\title{\textbf{Euclid's Number-Theoretical Work}}
\author{SHAOHUA ZHANG}
\date{{\small\emph{School of Mathematics, Shandong University,
Jinan,  Shandong, 250100, PRC\\
E-mail address: shaohuazhang@mail.sdu.edu.cn}}} \maketitle

\vspace{3mm}\textbf{Abstract:} When people mention the mathematical
achievements of Euclid, his geometrical achievements always spring
to mind. But, his Number-Theoretical achievements  (See Books 7, 8
and 9 in his magnum opus \emph{Elements} [1]) are rarely spoken. The
object of this paper is to affirm the number-theoretical role of
Euclid and the historical significance of Euclid's algorithm. It is
known that almost all elementary number-theoretical texts begin with
Division algorithm. However, Euclid did not do like this. He began
his number-theoretical work by introducing his algorithm. We were
quite surprised when we began to read the \emph{Elements} for the
first time. Nevertheless, one can prove that Euclid's algorithm is
essentially equivalent with the Bezout's equation and Division
algorithm. Therefore, Euclid has preliminarily established Theory of
Divisibility and the greatest common divisor. This is the foundation
of Number Theory. After more than 2000 years, by creatively
introducing the notion of congruence, Gauss published his
\emph{Disquisitiones Arithmeticae} in 1801 and developed Number
Theory as a systematic science. Note also that Euclid's algorithm
implies Euclid's first theorem (which is the heart of `the
uniqueness part' of the fundamental theorem of arithmetic) and
Euclid's second theorem (which states that there are infinitely many
primes). Thus, in the nature of things, Euclid's algorithm is the
most important number-theoretical work of Euclid. For this reason,
we further summarize briefly the influence of Euclid's algorithm.
Knuth said `we might call Euclid's method the granddaddy of all
algorithms'. Based on our discussion and analysis, it leads to the
conclusion  Euclid's algorithm is the greatest number-theoretical
achievement of the Euclidean age.

\vspace{3mm}\textbf{Keywords:}  Elements, Euclid's
number-theoretical work, Euclid's algorithm, Division algorithm,
Euclid's second theorem, Euclid's first theorem

\vspace{3mm}\textbf{2000 MR  Subject Classification:}00A05, 01A05,
11A41, 11A05


\section{ A brief introduction about Euclid's number-theoretical work}
\setcounter{section}{1}\setcounter{equation}{0}
Euclid, who was a Greek mathematician best known for his
\emph{Elements} which influenced the development of Western
mathematics for more than 2000 years, is one of the greatest
mathematicians of all time and popularly considered as the `Father
of Geometry', also known as Euclid of Alexandria, who lived probably
around 300 BC, is the most famous mathematician of antiquity. For
Euclid and the traditions about him, see [1, Introduction].

\vspace{3mm}His  magnum opus \emph{Elements} which covers much of
the geometry known to the ancient Greeks as well as some elementary
number theory (See \emph{Elements}: Books 7, 8 and 9) is probably
the most successful textbook ever written and has appeared in over a
thousand different editions from ancient to modern times. Heath [2,
Introduction] called it `…the greatest textbook of elementary
mathematics that there was written twenty-two centuries ago … Nor
does the reading of it require the 'higher mathematics'. Any
intelligent person with a fair recollection of school work in
elementary geometry would find it (progressing as it does by gradual
and nicely contrived steps) easy reading, and should feel a real
thrill in following its development, always assuming that enjoyment
of the book is not marred by any prospect of having to pass an
examination in it! …for everybody ought to read it who can, that is
all educated persons except the very few who are constitutionally
incapable of mathematics.' Heath [1, Preface] pointed out
`…Euclid's work will live long after all the text-books of the
present day are superseded and forgotten. It is one of the noblest
monuments of antiquity; no mathematician worthy of the name can
afford not to know Euclid, the real Euclid as distinct from any
revised or rewritten versions which will serve for schoolboys of
engineers. And, to know Euclid, it is necessary to know his
language, and so far as it can traced, the history of the 'elements'
which he collected in his immortal work.' We might not see such a
number-theoretical book any more, in which numbers are represented
by line segments and so have a geometrical appearance and aesthetic
feeling.

\vspace{3mm}From his \emph{Elements}, we know that Euclid's main
number-theoretical achievements should be reflected in the following
Propositions.

\vspace{3mm}\textbf{Proposition 1} (Book 7): Two unequal numbers
being set out, and the less being continually subtracted in turn
from the greater, if the number which is left never measures the one
before it until a unit is left, the original numbers will be prime
to one another.

\vspace{3mm}\textbf{Proposition 2} (Book 7): Given two numbers not
prime to one another, to find their greatest common measure.

\vspace{3mm}\textbf{Proposition 20} (Book 7): The least numbers of
those which have the same ratio with them measure those which have
the same ratio the same number of times, the greater the greater and
the lesser the lesser. Namely, if $a$ and $b$ are the smallest
numbers such that $a:b=c:d$, then $a$ divides $c$ and $b$ divides
$d$.

\vspace{3mm}\textbf{Proposition 30} (Book 7): If two numbers by
multiplying one another make some number, and any prime number
measures the product, it will also measure one of the original
numbers. Namely, if $p$ is prime, and $p|ab$, then $p|a$ or $p|b$.

\vspace{3mm}\textbf{Proposition 20} (Book 9): Prime numbers are more
than any assigned multitude of prime numbers. Namely, there are
infinitely many primes.

\vspace{3mm}\textbf{Proposition 36} (Book 9): If as many as we pleas
beginning from a unit be set out continuously in double proportion,
until the sum of all becomes prime, and if the sum multiplied into
the last make some number, the product will be perfect. Namely, if
$n=2^{p-1}(2^p-1)$, where $p$ is prime such that $2^p-1$  is also
prime, then, $n$ is even perfect number.

\vspace{3mm} Propositions 1 and 2 in Book 7 of \emph{Elements} are
exactly the famous Euclidean algorithm for computing the greatest
common divisor of two positive integers. According to Knuth [3], `we
might call Euclid's method the granddaddy of all algorithms, because
it is the oldest nontrivial algorithm that has survived to the
present day'.

\vspace{3mm}In their book \emph{An Introduction to the Theory of
Numbers}, Hardy and Wright [4] called  Proposition 20 (Book 9)
Euclid's second theorem. Hardy like particularly Euclid's proof of
for the infinitude of primes. Hardy [5] called it is `as fresh and
significant as when it was discovered---two thousand years have not
written a wrinkle on it'. According to Hardy [5], `Euclid's theorem
which states that the number of primes is infinite is vital for the
whole structure of arithmetic. The primes are the raw material out
of which we have to build arithmetic, and Euclid's theorem assures
us that we have plenty of material for the task'. Andr\'{e} Weil [6]
also called `the proof for the existence of infinitely many primes
represents undoubtedly a major advance, but there is no compelling
reason either for attributing it to Euclid or for dating them back
to earlier times. What matters for our purposes is that the very
broad diffusion of Euclid in later centuries, while driving out all
earlier texts, made them widely available to mathematicians from
then on'. I think that anyone who likes Number Theory must like
Euclid's second theorem. In his book The book of prime number
records, Paulo Ribenboim [7] cited nine and a half proofs of this
theorem. For other beautiful proofs, see [41--55].

\vspace{3mm}Hardy and Wright [4] called  Proposition 30 (Book 7)
Euclid's first theorem which is the heart of `the uniqueness part'
of the fundamental theorem of arithmetic. Recently, David Pengelley
and Fred Richman [8] published a readable paper entitled `Did Euclid
need the Euclidean algorithm to prove unique factorization'. They
called Proposition 30 (Book 7) Euclid's Lemma and pointed out that
Euclid's Lemma can be derived from Porism of Proposition 2. But `it
is not at all apparent that Euclid himself does this'. In their
paper, David Pengelley and Fred Richman explored that how Euclid
proved Proposition 30 using his algorithm. More precisely, they
proved Proposition 30 using the porism which states that if a number
divides two numbers, then it divides their greatest common divisor
as follows: if $p$ does not divide $a$, then  $\gcd (pb,ab)=b$, so
$p$ divides $b$. However, in their proof, they assumed such a clear
result $\gcd (pb,ab)=b$ which was not proven. Namely, in order to
prove that the porism of Proposition 2 implies that Proposition 30
holds, one need prove $\gcd (pb,ab)=b$ when $p$ does not divide $a$.
After giving Propositions 1, 2 and the porism of Proposition 2,
Euclid must want to prove that if   $\gcd (a,b)=1$ then  $\gcd
(ac,bc)=c$, which implies $\gcd (pb,ab)=b$ when $p$  does not divide
$a$. In fact, Euclid proved in spite of his expression is not like
this, that $\gcd (ac,bc)=c$ if and only if $\gcd (a,b)=1$, see [1]:
Book 7, Propositions 17, 18, 19, 20, 21(which states that numbers
prime to one another are the least of those which have the same
ratio with them) and 22 (which states that the least numbers of
those which have the same ratio with them are prime to one another).
So, Euclid's algorithm implies Euclid first theorem because his
proof of Proposition 30 refers exactly to Propositions 19 and 20.

\vspace{3mm}In most elementary number-theoretical texts, Euclid's
first theorem is derived from the Bezout's equation $ax+by=\gcd
(a,b)$. Of course, this is true (see also the following theorem 2).
Note that the Bezout's equation also can be derived from Euclid's
algorithm. It is a pity that Euclid himself does not obtain the
equation $ax+by=\gcd (a,b)$ by making use of his algorithm. If he
had known about $ax+by=\gcd (a,b)$! He perhaps knew, he just didn't
express explicitly.

\vspace{3mm}According to Hardy and Wright [4], the fundamental
theorem of arithmetic which says that every natural number is
uniquely a product of primes `does not seem to have been stated
explicitly before Guass. It was, of course, familiar to earlier
mathematicians; but Guass was the first to develop arithmetic as a
systematic science.' They further remarked: `It might seem strange
at first that Euclid, having gone so far, could not prove the
fundamental theorem itself; but this view would rest on a
misconception. Euclid had no formal calculus of multiplication and
exponentiation, and it would have been most difficult for him even
to state the theorem. He had not even a term for the product of more
than three factors. The omission of the fundamental theorem is in no
way casual or accidental; Euclid knew very well that the theory of
number turned upon his algorithm, and drew from it all the return he
could.' In his Disquisitiones Arithmeticae, Guass [9] proved
definitely that `a composite number can be resolved into prime
factors in only one way'. Maybe, Euclid did not consider this
problem how to prove unique factorization or how to resolve a
composite number into prime factors in only one way. But, he took
the first step and proved the following Propositions 31 and 32 by
the definition of composite number without using Euclid's algorithm
and its porism (see [1]: Book 7, the proofs of Propositions 31 and
32). Of course, he assumed that a finite composite number has only a
finite number of prime factors.

\vspace{3mm}\textbf{Proposition 31} (Book 7): Any composite number
is measured by some prime number.

\vspace{3mm}\textbf{Proposition 32 }(Book 7): Every number is either
prime or is measured by some prime number.

\vspace{3mm}About the problem on Euclid and the `fundamental theorem
of arithmetic', we have not pursued it. Knorr W. [56] gave a
reasonable discussion of the position of unique factorization in
Euclid's theory of numbers.

\vspace{3mm} Proposition 36 (Book 9) is Euclid's a great
number-theoretical achievement because he gave a sufficient
condition for even numbers to be perfect. In Weil's view, it is the
apex of Euclid's number-theoretical work [10]. A perfect number is
defined as a positive integer which is the sum of its proper
positive divisors, that is, the sum of the positive divisors
excluding the number itself. Euler proved further that the
sufficient condition about even perfect numbers given by Euclid 2000
years ago is also necessary. Namely, $n$ is even perfect number if
and only if $n=2^{p-1}(2^p-1)$, where $p$ is prime such that $2^p-1$
(also called Mersenne primes) is also prime. Thus, Euclid's work on
perfect numbers is not perfect. According to Littlewood [11],
`perfect number certainly never did any good, but then they never
did any particular harm.' Therefore, in this paper, we do not
further talk about Euclid's this work any more.

\vspace{3mm}Finally, we should mention again Proposition 20 in
Euclid's \emph{Elements} Book 7. B. L. van derWaerden [12] pointed
out that Proposition 20 plays a central role in Euclid's
arithmetical books. C. M. Taisbak [13] also announced Proposition 20
is the core of Euclid's arithmetical books. Proposition 20 can be
derived from Euclid's algorithm. Although `Central to Euclid's
development is the idea of four numbers being proportional:  $a$ is
to $b$ as $c$ is to $d$' [8], one can see again that the key of
studying divisibility is essentially Euclid's algorithm by David
Pengelley and Fred Richman's work. Generally speaking, Euclid's
algorithm is based on the following two results (Division algorithm
and Theorem 1) in Elementary Number Theory. As we know, Division
algorithm is the basis of Theory of Divisibility. Almost all
number-theoretical texts begin with it. For example, see [14~19].
However, Euclid did not do like this. Euclid began his
number-theoretical work by introducing his algorithm (See [1]: Book
7, Propositions 1 and 2). In this paper, we will prove  that
Euclid's algorithm is essentially equivalent with Division
algorithm. More precisely, for any positive integer $a$ and $b$,
that there exist unique integers $q$ and $r$ such that  $a=bq+r$ and
$0\leq r<b$ is equivalent with that there exist integers $x$ and $y$
such that $ax+by=\gcd (a,b)$. This implies that Division algorithm,
Euclid's algorithm and the Bezout's equation are equivalent. For the
details, see Section 2.

\section{Proof that Euclid's algorithm is equivalent with Division algorithm}
Strictly speaking, Division algorithm is essentially a theorem, but
Euclid's algorithm is an algorithm, a theorem and an algorithm are
not the same thing. Therefore, we should view Euclid's algorithm as
a theorem as the follows:

\vspace{3mm}\textbf{Euclid's algorithm:} For two distinct positive
integers, replace continually the larger number by the difference of
them until both are equal, then the answer is their greatest common
divisor.

\vspace{3mm}\textbf{Division algorithm (called also Division with
remainder):} For any positive integer $a$ and $b$, that there exist
unique integers $q$ and $r$ such that  $a=bq+r$ and $0\leq r<b$.

\vspace{3mm}\textbf{Proof of Division algorithm:} The uniqueness is
clear. Therefore, it is sufficient that we only prove the existence.
If $a|b$, then $q=\frac{b}{a}$ and $r=0$. Otherwise, consider the
set $S=\{a-bt:t=0,\pm1,\pm2,...,\infty\}$. It is manifest that there
are positive integers in this set. Therefore, there must be the
least positive integer in this set. Denote this number by
$c=a-bt_0$. Obviously, $c<b$(otherwise $c-b>0$ in the set $S$, and
$c$ is not the least positive integer in this set). Set $q=t_0$ and
$r=c$. This completes the proof.

\vspace{3mm}\textbf{Another proof of Division algorithm:} The
uniqueness is clear. Therefore, it is sufficient that we only prove
the existence. For any given positive integer $b$, when $a=1$, if
$b=1$, then let $q=1,r=0$, and the existence satisfies; if $b>1$,
then let $q=0,r=1$, and the existence satisfies again. Now, we
assume that the existence satisfies when $a=n$. We write
$n=bq_1+r_1$. Then, when $a=n+1$, we have $n+1=bq_1+r_1+1$. Note
that $0\leq r_1<b$. If $r_1=b-1$, then let $q=q_1+1,r=0$. If
$r_1<b-1$, then let $q=q_1,r=r_1+1$, and the existence satisfies
still. Therefore, by induction, Division algorithm is true.

\vspace{3mm}\noindent{\bf  Remark 1:~~}%
The proofs above need Peano axioms [14] which give the strict
definition of the set of natural numbers and imply the induction and
that if there are positive integers in a set, then, there must be a
least positive integer in this set and also imply the existence and
uniqueness the great common divisor.

\vspace{3mm}\textbf{Theorem 1:} For any positive integer $a$ and
$b$, that there exist unique integers $q$ and $r$ such that $a=bq+r$
and $0\leq r<b$. Moreover, if $a=bq+r$, then $\gcd (a,b)=\gcd
(b,r)$.

\vspace{3mm}Note that if  $b$ divides $a$, then $r=0$ in Division
algorithm and $\gcd (a,b)=b$. Hence, in order to accord with the
notation of the greatest common divisor of two positive integers, we
set $\gcd (b,r)=b$ when $r=0$.

\vspace{3mm}\textbf{ Proof of Theorem 1:} By Division algorithm, it
is sufficient that we only prove $\gcd (a,b)=\gcd (b,r)$. Let $\gcd
(a,b)=d$ and $\gcd (b,r)=t$. Since $r=a-bq$, and $d|a,d|b$, hence
$d|r$ . Now, $d|b$, $d|r$, so $d|\gcd (b,r)$ and $d|t$. On the other
hand, since $a=bq+r$, and $t|b$, $t|r$, hence $t|a$  and $t|d$. So,
$t=d$ and Theorem 1 is true.

\vspace{3mm}\noindent{\bf  Remark 2:~~}%
In this proof, we use the definition of the greatest common divisor
of positive integers $x$ and $y$, namely, the greatest common
divisor $\gcd (x,y)$ is not less than any other common divisor $x$
and $y$. Any two numbers have the greatest common divisor simply
because the set of common divisors is finite. But, this does not
imply that the greatest common divisor is divisible by any other
common divisor. By Peano axioms, there exists a least positive
integer $d$ in the set $S=\{ax+by:a\in Z,b\in Z\}$. By Division
algorithm, we have $x=dq+h$ and $y=dp+r$, but $r\in S,h\in S$, so
$r=h=0$ since $d$ is least. Thus $d$ is a common divisor $x$ and
$y$. But $d=ax+by$ implies that $d$ is divisible by any other common
divisor of $x$ and $y$. Specially,  $d$ is not less than any other
common divisor of $x$ and $y$. By the definition of the greatest
common divisor, $d$ is the greatest common divisor. The uniqueness
of the greatest common divisor is clear. It immediately shows that
the greatest common divisor of $x$ and $y$ is divisible by any other
common divisor of $x$ and $y$. This is exactly the porism of
Proposition 2 in \emph{Elements} Book 7.

\vspace{3mm}Nevertheless, in his Elements Book 7, Euclid did not
give the definition of the greatest common divisor. He only gave the
definition of numbers prime to one another, in this case, the
greatest common divisor of numbers is 1, see [1: Book 7, definition
12]. Then, he gave the sufficient condition of two numbers prime to
one another, see [1: Book7, Proposition 1]. I wonder why he did not
give a necessary condition or the necessary and sufficient condition
of two numbers prime to one another as Proposition 1. More
precisely, Proposition 1 should be expressed as follows: two numbers
will be prime to one another if and only if the less being
continually subtracted in turn from the greater, if the number which
is left never measures the one before it until an unite (i.e.1) is
left. After proving Proposition 1, Euclid considered how to find the
greatest common divisor of two numbers, see [1: Book 7, Proposition
2]. He proved firstly that the output of his algorithm is a common
divisor of two numbers, and then, he showed that any common divisor
has to divide it, so has to be smaller. Thus, he gave naturally the
porism of Proposition 2 which states that if a number is a common
divisor of two numbers, then it divides their greatest common
divisor.

\vspace{3mm}By Division algorithm and Theorem 1, we see again that
Euclid's algorithm is true without using Euclid's proof in
Proposition 2. So, we can say that Division algorithm implies
Euclid's algorithm by giving the definition of the greatest common
divisor. In fact, from his algorithm (Propositions 1 and 2), one can
observe that Euclid assumed as two axioms. One is of that, if $a$
and  $b$ are both divisible by $c$, so is $a-bq$. Another is of that
for any positive integer $a$ and $b$, there exist integers $q$ and
$r$ such that $a=bq+r$ and $0\leq r<b$. We do not know why Euclid
began his number-theoretical work by introducing his algorithm
without introducing any postulates or common notions (axioms). This
is not his style for constructing Geometry system. He must have some
reason for this. In my eyes, Euclid is not only a mathematician but
also a philosopher. Perhaps, he felt that his algorithm is
essentially equivalent with Division algorithm, or that he needed
more complicated expression for some axioms like Peano axioms.

\vspace{3mm}Now, we consider the set $T=\{ax+by:x\in Z,y\in Z\}$ for
any given positive integers $a$ and $b$. Note that there are
positive integers in this set. Therefore, there must be a least
positive integer in this set. Denote this number by $d=ax_0+by_0$.
Of course, $d$  is just the greatest common divisor of $a$ and $b$.
Thus we get the following Bezout's equation.

\vspace{3mm}\textbf{Theorem 2: }For any given positive integers $a$
and $b$, there are integers $x$ and $y$ such that $ax+by=\gcd
(a,b)$.

\vspace{3mm} By Euclid's algorithm, it is also easy to find such
integers $x$ and  $y$ so that $ax+by=\gcd (a,b)$. In fact, we have
the well-known procedure involves repeated division, resulting in a
sequence as following: $a=bq_1+r_1$ ($0\leq r_1<b$), $b=q_2r_1+r_2$
($0\leq r_2<r_1$), $r_1=q_3r_2+r_3$ ($0\leq r_3<r_2$),...., The
process terminates when some remainder $r_k=0$, which must happen
eventually. Then $\gcd (a,b)=r_{k-1}$, and $x$ and $y$ are found as
follows: Compute $r_1=a-q_1b$, substitute into $b=q_2r_1+r_2$ and
obtain $r_2=-q_2a+(1+q_1q_2b)$, and so on. Eventually, we must find
integers $x$ and $y$ such that $ax+by=\gcd (a,b)$. So, Euclid's
algorithm implies Theorem 2.

\vspace{3mm} Next, we prove that Theorem 2 implies Division
algorithm. By symmetry it is safe to assume that $x\geq 0$. If
$x=0$, then $b|a$, so there exist integers $q$ and $r$ such that
$a=bq+r$ and $0= r<b$. If $x=1$, then $a=-yb+\gcd (a,b)$. If $\gcd
(a,b)=b$, then $b|a$. So we can assume that $\gcd (a,b)<b$. Let
$q=-y, r=\gcd (a,b)$, then Division algorithm holds. Now let $x>1$.
Since $\gcd (a,b)=ax+by$, hence we can write
$a=b(-y-x+1)+(x-1)(b-a)+\gcd (a,b)$. If $0\leq (x-1)(b-a)+\gcd
(a,b)<b$, then let $r=(x-1)(b-a)+\gcd (a,b)$ and Division algorithm
holds. If $(x-1)(b-a)+\gcd (a,b)> b$, then $b>a$ since $x>1$ and
$\gcd (a,b)\leq b$. Let $q=0$ and $r=a$. Division algorithm holds.
If $(x-1)(b-a)+\gcd (a,b)= b$, then $b|a$ since
$a=b(-y-x+1)+(x-1)(b-a)+\gcd (a,b)$.  Division algorithm holds
still. Therefore, it suffices to consider the case of
$(x-1)(b-a)+\gcd (a,b)<0$ which implies that $0<b<a-\frac{\gcd
(a,b)}{x-1}<a$ . We write $(x-1)(b-a)+\gcd (a,b)=-c$. Thus,
$0<c<b(-y-x+1)$. Let $d=-y-x+1$. Clearly, $d$ is a positive integer.
Note that $x$ and $y$ are given since $ax+by=\gcd (a,b)$. Therefore
$d$ is decided. So, we can consider a finite number of intervals
$[0,b),[b,2b),...,[(d-1)b,db)$. $c$ must be in some interval among
these intervals. Let $ib\leq c<(i+1)b$, where $1\leq i<d$. Thus,
$a=bd-c=b(d-i-1)+(i+1)b-c$. Let $q=d-i-1$ and $r=(i+1)b-c$.
Therefore, there exist integers $q$ and $r$ such that $a=bq+r$ and
$0\leq r<b$. Of course, the uniqueness of $q$ and $r$ is obvious.
Thus, Division algorithm is true. Therefore, Euclid's algorithm is
essentially equivalent with Division algorithm. Thus, we proved the
following theorem 3.

\vspace{3mm}\textbf{Theorem 3: }For any positive integer $a$ and
$b$, that there exist unique integers $q$ and $r$ such that $a=bq+r$
and $0\leq r<b$ is equivalent with that there exist integers $x$ and
$y$ such that $ax+by=\gcd (a,b)$.

\vspace{3mm}Theorem 3 shows that Division algorithm, Euclid's
algorithm and the Bezout's equationSee are equivalent.

\vspace{3mm}Note that  Euclid's algorithm implies also Euclid's
second theorem. Let's go back to Euclid's proof for the infinitude
of prime numbers: Supposed that there are only finitely many primes,
say $k$ of them, which denoted by $p_1,...,p_k$. Consider the number
$E=1+ \prod_{i=1}^{i=k}p_i$. If $E$ is prime, it leads to the
contradiction since $E\neq p_i$ for any $1\leq i\leq k$. If $E$ is
not prime, $E$ has a prime divisor $p$ by Proposition 31 (Book 7).
But $p\neq p_i$ for any $1\leq i\leq k$. Otherwise, $p$ divides
$\prod_{i=1}^{i=k}p_i$. Since it also divides
$1+\prod_{i=1}^{i=k}p_i$, it will divide the difference or unity,
which is impossible.

\vspace{3mm} In his proof, we see that Euclid used Proposition 31
(Book 7). Of course, he also used a unexpressed axiom which states
that if $A$ divides $B$, and also divides $C$, $A$ will divide the
difference between $B$ and $C$.

\vspace{3mm}Well, let's look at the proof of Proposition 31 (Book
7): Let $A$ be a composite number. By the definition, there must be
a number $B$ ($1<B<A$) which divides $A$. If $B$ is prime, then
Proposition 31 holds. If  $B$ is not prime, there must be a number
$C$ ($1<C<B$) which divides $B$. If $C$ is prime, then Proposition
31 holds since $C$ also divides $A$. If $C$ is not prime, by
repeating this process, in  finite many steps, there must be a prime
which divides $A$ and Proposition 31 holds. From this proof, we see
that Euclid used a unexpressed axiom which states that if $A$
divides $B$, and $B$ divides $C$, then $A$ divides $C$. Thomas
Little Heath had noted that Euclid used the aforementioned axioms.
We would be quite surprised if he did use these axioms because  on
one hand, Proposition 31 (Book 7) and Proposition 20 (Book 9) can be
deduced early by definitions, on the other hand, we expect him to
make use of his algorithm which is his first number-theoretical
proposition in his \emph{Elements}.  In [69, Appendix], we try to
supplement this work and give a simple proof.

\vspace{3mm} For his other number-theoretical work, see [Appendix]
in which we give a list of 22 definitions, 102 propositions and 3
porisms on Number Theory in Euclid's  \emph{Elements} Books 7, 8 and
9 (Translated by Thomas Little Heath). We omit Euclid's proofs and
Heath's commentaries for anyone who would like to learn Euclid's
Number Theory in his own way.

\vspace{3mm}Based on the discussion and analysis above, we can say
that we ourselves `require the Euclidean algorithm to prove unique
factorization' because it is the basis of Theory of Divisibility
from Peano axioms. \emph{Elements} established elementarily Theory
of Divisibility and the greatest common divisor. This is the base of
Number Theory (After more than 2000 years, by introducing creatively
the notion of congruence, Gauss published his \emph{Disquisitiones
Arithmeticae} in 1801 and developed Number Theory as a systematic
science.). In the nature of things, Euclid's algorithm is the most
important number-theoretical achievement of Euclid. In next section,
we will further summarize briefly the influence of Euclid's
algorithm.

\section{Euclid's Algorithm and Our Conclusions}
Historically, many mathematicians and computer scientists studied
Euclid's algorithm. For example, D. H. Lehmer, J. D. Dixon, L. K.
Hua, Donald E. Knuth, Andrew C. Yao, H. Lenstra, A. K. Lenstra, H.
Davenport, J. Barkley Rosser, P.M. Cohn, Heilbronn, Viggo Brun and
so on. So many people like Euclid's because it is not only simple
and beautiful but also useful. Although more than 2000 years have
passed, the study on Euclid's algorithm still goes on and on. Heath
[1, Preface to the second edition] said: `…Euclid is far from being
defunct or even dormant, and that, so long as mathematics is
studied, mathematicians will find it necessary and worth while to
come back again and again, for one purpose or another, to the
twenty-two centuries-old book which, notwithstanding its
imperfections, remains the greatest elementary textbook in
mathematics that the world is privileged to possess.'

\vspace{3mm}In 1968, Heilbronn [22] studied the average length of a
class of finite continued fractions. This is an important result on
Euclid’s algorithm. T. Tonkov [23] and J. W. Porter [24] improved
Heilbronn’s estimate respectively. In 1975, using an idea of
Heilbronn, Andrew C. Yao and Donald E. Knuth [25] studied the sum of
the partial quotients $q_i$ in Euclid’s algorithm. They proved a
well-known result which states that the sum $S$ of all the partial
quotients of all the regular continued fractions for $\frac{a}{b}$
with $1\leq b\leq a$ is $6\pi ^{-2}a (\log a)^2+ O(a\log a( \log\log
a)^2)$. This implies that $S\ll a(\log a)^2$. In 1994, ZhiYong Zheng
[26] improved the result of Andrew C. Yao and Donald E. Knuth. As an
application, in 1996, J. B. Conrey, Eric Fransen, and Robert Klein
[27] studied the mean values of Dedekind Sums. For a positive
integer $k$ and an arbitrary integer $h$, the Dedekind sum $s(h,k)$
is defined by
$s(h,k)=\sum_{a=1}^{a=k}((\frac{a}{k}))((\frac{ah}{k}))$, where
$((x))=x-[x]-\frac{1}{2}$ if $x\neq 0$, otherwise $((x))=0$. The
most famous property of the Dedekind sums is the reciprocity formula
as follows $s(h,k)+s(k,h)=\frac{h^2+k^2+1}{12hk}-\frac{1}{4}$.
Dedekind Sums are closely related to the transformation theory of
the Dedekind eta-function which is the infinite product $\eta
(\tau)=q_{24}\prod _{n=1}^{n=\infty}(1-q^n)$ and satisfies the
transformation law $\eta (\frac{-1}{\tau})=\sqrt{-i\tau}\eta (\tau)$
with $\tau\in H= \{\tau\in C:Im (\tau)>0\}$,  where $q=e^{2\pi
i\tau}, q_{24}=e^{\frac{1}{12}\pi i\tau}$. Thus, the Dedekind
eta-function is a modular form of weight $\frac{1}{2}$ for the
modular group $SL(2,Z)=\{\left(
\begin{matrix}
a, b\\
c,d
\end{matrix} \right):a,b,c,d\in Z,ad-bc=1\}$. Each element of the modular group $SL(2,Z)$ is
also viewed as an automorphism (invertible self-map) of the Riemann
sphere $C\cup\{\infty\}$, the fractional linear transformation
$\left(
\begin{matrix}
a, b\\
c,d
\end{matrix} \right)(\tau)=\frac{a\tau+b}{c\tau+d},
\tau\in C\cup\{\infty\}$.     The Dedekind eta-function played a
prominent role in Number Theory and in other areas of mathematics.
From these results above, one can see again the importance of
Euclid's algorithm.

\vspace{3mm}    Note that Euclid’s algorithm is essentially a
dynamical system. Namely, the Euclidean algorithm is the map defined
by $$f:(x,y)\in R_+^2\rightarrow \left\{
\begin{array}{c}
(x-y,y), x\geq y
\\
(x,y-x),x< y
\end{array}
\right.,$$ where $R_+=\{x\in R:x\geq 0\}$. Clearly, when $x$ and $y$
are natural numbers, there must be a positive integer $k$ such that
$f^k(x,y)=(0,d)$ or $(d,0)$, where $d$ is the greatest common
divisor of $x$ and $y$. Equivalently, the map $f$ is given by $$
\left(
\begin{matrix}
x\\
y
\end{matrix} \right)\rightarrow \left\{
\begin{array}{c}
\left(
\begin{matrix}
1, -1\\
0,1
\end{matrix} \right)\left(
\begin{matrix}
x\\
y
\end{matrix} \right), x\geq y
\\
\left(
\begin{matrix}
1, 0\\
-1,1
\end{matrix} \right)\left(
\begin{matrix}
x\\
y
\end{matrix} \right),x< y
\end{array}
\right..$$

Note that $\left(
\begin{matrix}
1, -1\\
0,1
\end{matrix} \right)^{-1}$ and $\left(
\begin{matrix}
1, 0\\
-1,1
\end{matrix}\right)^{-1}$ generate the modular group $SL(2,Z)$.     By a result of Hedlund [28],
we know that the linear action of $SL(2,Z)$ on $R^2$ is ergodic. By
a result of A. Nogueira [29], we see also that the map $f$ is
ergodic relative to the Lebesgue measure on $R^2_+$. In 2007, Dani,
S. G. and Nogueira, Arnaldo [30] showed that `the Euclidean
algorithm $f$ turns out to be an example of a dissipative
transformation for which this is not the case via a natural
extension of $f$ constructed by using the action of $SL(2,Z)$ on a
subset of $SL(2,R)$'.

 \vspace{3mm}   Euclid's algorithm gives the greatest common divisor $a$ and $b$ in a finite number of steps
 and its total running time is approximately $O(\log ab)$ word operations [31]. Lam\'{e}'s theorem states
 that Euclid's algorithm runs in no more than $5k$ steps, where $k$ is the number of (decimal) digits of $\max \{a,b\}$ [17].
 In 2003, Vall\'{e}e, B. [32] published a paper entitled `Dynamical analysis of a class of Euclidean algorithms' and
 developed a general technique for analyzing the average-case behavior of the Euclidean-type algorithms. `This is a deep and
 important paper which merits careful study, and will likely have a significant impact on future directions in algorithm analysis.'
  Jeffrey O. Shallit reviewed, `The method involves viewing these algorithms as a dynamical system, where each step is a linear
  fractional transformation of the previous one. ... Then a generating function (Dirichlet series) is used to describe the cost of
  the algorithms, and Tauberian theorems are used to extract the coefficients.' In 2006, Vall\'{e}e, B. [33] further
  proposed a detailed and precise dynamical and probabilistic analysis of the more natural variants of Euclid's algorithm.
  The paper `presented a clear, clever, and unified overview of the methodology (dynamical analysis) and of the tools.'
  Val\'{e}rie Berth\'{e} reviewed, `... One of the strengths of this approach comes from the fact that it combines
  sophisticated tools taken from, on the one hand, analytic combinatorics and functional analysis
  (moment generating functions, Dirichlet series, quasi-power theorems and Tauberian theorems),
  and on the other hand, from dynamical systems and ergodic theory (including Markovian dynamical
  systems, induction, transfer operators and related concepts such as zeta series, pressure and entropy)...
  In particular, Section 9 provides a detailed overview of the literature on the subject as well
  as a list of open problems, particularly, in the higher-dimensional case (including the lattice
  reduction problem or classical multidimensional continued fraction algorithms such as the Jacobi-Perron algorithm).'
  From these results above, one would see again the vitality of Euclid's algorithm. `Although more than 2000 years have passed,
  the Euclidean algorithm has not yet been completely analyzed'---Vall\'{e}e, B.

 \vspace{3mm}     In Book 7, Proposition 3 of his Elements [1],
 Euclid considered how to compute the great common divisors of three positive integers $a,b$ and $c$.
 His method is simple and natural. Namely, firstly, compute $\gcd (a,b)=d$, secondly, compute $\gcd(d,c)=e$, then $\gcd (a,b,c)=e$.
 It is easy to prove that this method is true, and this method can be readily generalized to the case for computing the greatest
 common divisor of several positive integers. Thus, using Euclid's algorithm, one can solve the following two problems:

 \vspace{3mm}\textbf{Problem 1:} Let $a_1,...,a_m$ and $b$ be any positive integers.
 Find an algorithm to determine whether $b$ can be represented by $a_1,...,a_m$. Or equivalently, find an algorithm to determine whether $b$
 belongs to the ideal generated by $a_1,...,a_m$.

 \vspace{3mm}\textbf{Problem 2:} Let $a_1,...,a_m$ be any positive integers.
 Find an algorithm to determine whether there is an integer $a_i$ among  $a_1,...,a_m$ such that $a_i$ is relatively prime with all of the others.

  \vspace{3mm}    Generalizing Problem 1 to the case over the unique factorization domain $F[x_1,...,x_n]$, where $F$ is a field.
  Let $f_1,...,f_m$ and $g$ be polynomials in $F[x_1,...,x_n]$. Find an algorithm to determine whether $g$ belongs to the ideal generated by $f_1,...,f_m$.
   As we know, this interesting generalization leads to the invention of Gr\"{o}bner bases of polynomial ideals.
   Using S-polynomial, also combining with the multivariate division algorithm, in 1965, Buchberger [34]
   gave an algorithm for finding a basis $g_1,...,g_k$ of the ideal $I=<f_1,...,f_m>$ such that the leading term
   of any polynomial in $I$ is divisible by the leading term of some polynomial in $G=\{g_1,...,g_k\}$.
   Such a basis is called Gr\"{o}bner bases by Buchberger. An analogous concept was developed independently
   by Heisuke Hironaka in 1964 [35, 36], who named it standard bases.
   As we know, Gr\"{o}bner bases of polynomial ideals in modern Computational Algebraic
   Geometry are very important. They also are rather useful to Symbolic Computation and Cryptography and so on.
   The concept and algorithms of Gr\"{o}bner bases have been further generalized to modules
   over a polynomial ring and to non-commutative (skew) polynomial rings such as Weyl algebras.

  \vspace{3mm}        Problem 2 leads to the invention of W sequences [58] which play an interesting role in the study
  of primes and enable us to give new weakened forms of many classical problems which are open in Number Theory.
  For any integer $n>1$, the sequence of integers $0<a_1<a_2<...<a_n$ is called a W sequence, if there exists $r$ with $1\leq r\leq n$ such that
  $a_r$ and each of the rest numbers are coprime.

  \vspace{3mm}    W sequences in the case of consecutive positive integers relates to Grimm's conjecture.
  Clearly, if there is a prime in the sequence $m+1,...,m+n$, then this sequence is a W sequence.
  Therefore, in order to determine whether a consecutive positive integer sequence is a W sequence,
  it is enough to consider the case of consecutive composite numbers.
  It leads to the further study of consecutive composite numbers.
  In 1969, C.A.Grimm [57] made an important conjecture that if  $m+1,...,m+n$ are consecutive composite numbers,
  then there exist $n$ distinct prime numbers $p_1,...,p_n$ such that $m+i$ is divisible by $p_i$ for $1\leq i\leq n$.
   This implies that for all sufficiently large integer $n$, there is a prime between $n^2$ and $(n+1)^2$.
   It is nice that for $m\geq 1$, that there exists a prime in the interval $(m^2,(m+1)^2)$ is equivalent with that $m^2+1,...,m^2+2m$
   is a W sequence. In [59], we further refine the function $g(m)$ on Grimm's conjecture and obtain several interesting results.
   For example, we refine a result of Erd\"{o}s and Selfridge without using Hall's theorem.

\vspace{3mm}  Denote the largest integer $n$  in  $m+1,...,m+n$ by
$h(m)$ such that no one of  $m+1,...,m+n$ is relatively prime with
all of the others. Cram\'{e}r's conjecture [70] and Pillai's result
[71] imply $17\leq h(m)\ll ( \log m)^2$. From these, one can see
that the W sequences in the case of consecutive positive integers
tie up the distribution of primes in short interval. Unfortunately,
it is not easy to prove that a sequence is a W sequence.

\vspace{3mm}In the non-consecutive case, W sequences enable us to
get a new weakened form of Goldbach's conjecture and reveal the
internal relationship between Goldbach's conjecture and the least
prime in an arithmetic progression. We find that Goldbach Conjecture
ties up Kanold Hypothesis and Chowla Hypothesis, for details, see
[60]. For another example, for positive integers $a$ and $b$ with
$\gcd (a,b)=1$ and $1\leq a<b$, if there is a prime in
$a,a+b,...,a+(b-1)b$, then this sequence is a W sequence. Thus, in
order to determine whether the sequence $a,a+b,...,a+(b-1)b$ is a W
sequence, it is enough to consider the case that
$a,a+b,...,a+(b-1)b$ are all composite numbers. It leads to the
generalization of Grimm's conjecture. In [61], we generalized a
theorem about the binomial coefficient and got a slightly stronger
result than Langevin's [62]. This leads to possible generalizations
of Grimm's conjecture [68].

\vspace{3mm} Generalize W sequences to the case over the unique
factorization domain $F[x_1,...,x_n]$, where $F$ is a field. It
leads how to find an efficient algorithm for computing the greatest
common divisor of any two polynomials in $F[x_1,...,x_n]$.

\vspace{3mm} By the aforementioned discussion, one see that Euclid's
algorithm implies his second theorem. We believe that one of
substantive characteristics of the set of all integers is that it
contains infinitely many prime numbers. It is known that $f(x)=x$ on
$Z$ is the simplest polynomial which represents infinitely many
primes. By  Dirichlet's famous theorem, for any positive integer
$l,k$ with $(l,k)=1$, $f(x)=l+kx$ is a simpler polynomial which also
represents infinitely many primes. More generally, it is possible to
establish a generic model for the problem that several multivariable
integral polynomials represent simultaneously primes. More
concretely, let's consider the map $F: Z^n\rightarrow Z^m$ for all
integral points $x=(x_1,...,x_n)\in Z^n$, $F(x)=(f_1(x),...,f_m(x))$
for distinct polynomials $f_1,...,f_m\in Z[x_1,...,x_n]$. How to
determine whether $f_1(x),...,f_m(x)$ represent simultaneously
primes? In [67,72], we considered this problem and obtained some
interesting results. We strongly wish that in the higher-dimension
case, we have a similar Prime Number Theorem. Thus, it is also
possible to generalize the problem of the least prime number in an
arithmetic progression and  an analogy of Chinese Remainder Theorem,
moreover, give an analogy of Goldbach's conjecture, and so on.

\vspace{3mm}It is well-known that we can solve either linear
Diophantine equations or a system of simultaneous linear
congruences, and find also modular inversions, and expand continued
fractions, testing primality, generating primes, factoring large
integers and so on by using Euclid's algorithm and extended
algorithm. Moreover, due to the fact that Euclid's algorithm can be
not only readily generalized to polynomials in one variable over a
field but also generalized multivariate polynomials over any unique
factorization domain, Euclid's algorithm also plays an important
role in symbolic computation and cryptography even in science and
engineering. Without Euclid's algorithm there would be no the
prosperity of computation nowadays, we are afraid. It is very nice
that Viggo Brun [20, 21] studied the relations between Euclid's
algorithm and music theory. By coincidence, Euclid himself also
reveled in music [1].

\vspace{3mm}Looking back into Ancient Greek Number Theory history,
it is not difficult to confirm that Euclid's algorithm indeed is the
greatest number-theoretical achievements of the age. Some scholars
believe that Euclid's algorithm probably was not Euclid's own
invention [19]. Many scholars conjecture was actually Euclid's
rendition of an algorithm due to Eudoxus (c.375B.C.) [3].
Nevertheless, it first appeared in Euclid's Elements, and more
importantly, it is the first nontrivial algorithm that has survived
to this day. Therefore, I think that this is why people would like
to call it Euclid's algorithm. Perhaps, it is suitable to call it
Ancient Greek Algorithm.

\vspace{3mm}Anyway, also closely relating to many famous algorithms
such as Guass' elimination, Buchberger's algorithm [34], Schoof's
algorithm [37, 38], Cornacchia's algorithm [39], LLL Algorithm [40],
modern factorization algorithms (Continued Fraction Factoring
Algorithm [63], the Elliptic Curve Factoring Algorithm [64], the
Multiple Polynomial Quadratic Sieve [65] and the Number Field Sieve
[66]) and so on, the Euclidean algorithm (together with the
discovery of irrationals in Pythagoras' School) is the greatest
achievement of Ancient Greek Number Theory.  Let's cite Edna St.
Vincent Millay' sonnet Euclid Alone Has Looked on Beauty Bare to
close this paper.

\vspace{3mm}\textbf{Euclid alone has looked on Beauty bare}

\vspace{3mm}Euclid alone has looked on Beauty bare.

Let all who prate of Beauty hold their peace,

And lay them prone upon the earth and cease

To ponder on themselves, the while they stare

At nothing, intricately drawn nowhere

In shapes of shifting lineage; let geese

Gabble and hiss, but heroes seek release

From dusty bondage into luminous air.

O blinding hour, O holy, terrible day,

When first the shaft into his vision shone

Of light anatomized! Euclid alone

Has looked on Beauty bare. Fortunate they

Who, though once only and then but far away,

Have heard her massive sandal set on stone.

\section*{Acknowledgements}
I am happy to record my gratitude to my supervisor Professor Xiaoyun
Wang who made a lot of helpful suggestions. I wish to thank the
referees for a careful reading of this paper. Thank Benno van Dalen
for his valuable reviews. Thank Professors Guangwu Xu, Mingqiang
Wang for their help. Thank the Institute for Advanced Study in
Tsinghua University for providing me with excellent conditions. This
work was partially supported by the National Basic Research Program
(973) of China (No. 2007CB807902) and the Natural Science Foundation
of Shandong Province (No.Y2008G23). Y2008G23).

\section*{Appendix: Euclid's \emph{Elements} Books 7, 8 and 9}

Below is the list of 22 definitions, 102 propositions and 3 porisms
on Elementary Number Theory in Euclid's \emph{Elements} Books 7, 8
and 9 (Translated by Thomas Little Heath).

\vspace{3mm}\textbf{Book 7 Definitions }

1. A unit is that by virtue of which each of the things that exist
is called one.

2. A number is a multitude composed of units.

3. A number is a part of a number, the lesser of the greater, when
it measures the greater.

4. but parts when it does not measure it.

5. The greater number is a multiple of the lesser when it is
measured by the lesser.

6. An even number is that which is divisible into two equal parts.

7. An odd number is that which is not divisible into two equal parts
or that which differs by a unit from an even number.

8. An even-times-even number is that which is measured by an even
number according to an even number.

9. An even-times-odd number is that which is measured by an even
number according to an odd number.

10. An odd-times-odd number is that which is measured by an odd
number according to an odd number.

11. A prime number is that which is measured by a unit alone.

12. Numbers prime to one another are those which are measured by a
unit alone as a common measure.

13. A composite number is that which is measured by some number.

14. Numbers composite to one another are those which are measured by
some number as a common measure.

15. A number is said to multiply a number when that which is
multiplied is added to itself as many times as there are units in
the other, and thus some number is produced.

16. And, when two numbers having multiplied one another make some
number, the number so produced is called plane, and its sides are
the numbers which have multiplied one another.

17. And, when three numbers having multiplied one another make some
number, the number so produced is solid, and its sides are the
numbers which have multiplied one another.

18. A square number is equal multiplied by equal, or a number which
is contained by two equal numbers.

19. And a cube is equal multiplied by equal and again by equal, or a
number which is contained by three equal numbers.

20. Numbers are proportional when the first is the same multiple, or
the same part, or the same parts, of the second that the third is of
the fourth.

21. Similar plane and solid numbers are those which have their sides
proportional.

22. A perfect number is that which is equal to its own parts.

\vspace{3mm}\textbf{Book 7 Propositions }

Proposition 1: Two unequal numbers being set out, and the less being
continually subtracted in turn from the greater, if the number which
is left never measures the one before it until a unit is left, the
original numbers will be prime to one another.

Proposition 2: Given two numbers not prime to one another, to find
their greatest common measure.

Porism: If a number measures two numbers then it will also measure
their greatest common measure.

Proposition 3: Given three numbers not prime to one another, to find
their greatest common measure.

Proposition 4: Any number is either a part or parts of any number,
the lesser of the greater.

Proposition 5: If a number be a part of a number, and another be the
same part of another, then the sum will also be the same part of the
sum that the one is of the one.

Proposition 6: If a number be parts of a number, and another be the
same parts of another, then the sum will also be the same parts of
the sum that the one is of the one.

Proposition 7: If a number be that part of a number, which a number
subtracted is of a number subtracted, the remainder will also be the
same part of the remainder that the whole is of the whole.

Proposition 8: If a number be the same parts of a number that a
number subtracted is of a number subtracted, the remainder will also
be the same parts of the remainder that the whole is of the whole.

Proposition 9: If a number be a part of a number, and another be the
same part of another, alternately also, whatever part or parts the
first is of the third, the same part, or the same parts, will the
second also be of the fourth.

Proposition 10: If a number is parts of a number, and another be the
same parts of another, alternately also, whatever parts or part the
first is of the third, the same parts, or the same part will the
second also be of the fourth.

Proposition 11: If, as the whole is to whole, so is a number
subtracted to a number subtracted, the remainder will also be to the
remainder as whole to whole.

Proposition 12: If there be as many numbers as we please in
proportion, then, as one of the antecedents is to one of the
consequents, so are all the antecedents to all the consequents.

Proposition 13: If four numbers be proportional, they will also be
proportional alternately.

Proposition 14: If there be as many numbers as we please, and others
equal to them in multitude, which taken two and two are in the same
ratio, they will also be in the same ratio ex aequali.

Proposition 15: If a unit measures any number, and another number
measures any other number the same number of times, then,
alternately also, the unit will measure the third number the same
number of times that the second measures the fourth.

Proposition 16: If two numbers by multiplying one another make
certain numbers, the numbers so produced will be equal to one
another.

Proposition 17: If a number by multiplying two numbers make certain
numbers, the numbers so produced will have the same ratio as the
number multiplied.

Proposition 18: If two numbers by multiplying any number make
certain numbers, the numbers so produced will have the same ratio as
the multipliers.

Proposition 19: If four numbers be proportional, the number produced
from the first and fourth will be equal to the number produced from
the second and third; and, if the number produced from the first and
fourth be equal to that produced from the second and third, the four
numbers will be proportional.

Proposition 20: The least numbers of those which have the same ratio
with them measure those which have the same ratio the same number of
times, the greater the greater and the lesser the lesser.

Proposition 21: Numbers prime to one another are the least of those
which have the same ratio with them.

Proposition 22: The least numbers of those which have the same ratio
with them are prime to one another.

Proposition 23: If two numbers be prime to one another, the number
which measures the one of them will be prime to the remaining
number.

Proposition 24: If two numbers be prime to any number, their
produced also will be prime to the same.

Proposition 25: If two numbers be prime to one another, the product
of one of them into itself will be prime to the remaining one.

Proposition 26: If two numbers be prime to two numbers, both to
each, their products also will be prime to one another.

Proposition 27: If two numbers be prime to one another, and each by
multiplying itself make a certain number, the products will be prime
to one another; and if the original numbers by multiplying the
products make certain numbers, the latter will also be prime to one
another [and this is always the case with the extremes].

Proposition 28: If two numbers be prime to one another, the sum will
also be prime to each of them; and if the sum of two numbers be
primes to any one of them, the original numbers will also be prime
to one another.

Proposition 29: Any prime number is prime to any number which it
does not measure.

Proposition 30: If two numbers by multiplying one another make some
number, and any prime number measures the product, it will also
measure one of the original numbers.

Proposition 31: Any composite number is measured by some prime
number.

Proposition 32: Every number is either prime or is measured by some
prime number.

Proposition 33: Given as many numbers as we please, to find the
least of those which have the same ratio with them.

Proposition 34: Given two numbers, to find the least number which
they measure.

Proposition 35: If two numbers measure any number, the least number
measured by them will also measure the same number.

Proposition 36: Given three numbers, to find the least number which
they measure.

Proposition 37: If a number be measured by any number, the number
which is measured will have a part called by the same name as the
measuring number.

Proposition 38: If a number has any part whatever, it will be
measured by a number called the same name as the part.

Proposition 39: To find the number which is the least that will have
given parts.

\vspace{3mm}\textbf{Book 8 Propositions }

Proposition 1: If there be as many numbers as we please in continued
proportion, and the extremes of them be prime to one another, the
numbers are the least of those which have the same ratio with them.

Proposition 2: To find the numbers in continued proportion, as many
as may be prescribed, and the least that are in a given ratio.

Porism: If three numbers in continued proportion be the least of
those which have the same ratio with them, the extremes of them are
square, and, if four numbers, cubes.

Proposition 3: If as many numbers as we please in continued
proportion be the least of those which have the same ratio with
them, the extremes of them are prime to one another.

Proposition 4: Given as many numbers as we please in least numbers,
to find numbers in continued proportion which are the least in the
given ratios.

Proposition 5: Plane numbers have to one another the ratio
compounded of the ratios of their sides.

Proposition 6: If there be many numbers as we please in continued
proportion, and the first does not measure the second, then neither
will any other measure any other.

Proposition 7: If there be as many numbers as we please in continued
proportion, and the first measures the last, it will measure the
second also.

Proposition 8: If between two numbers there fall numbers in
continued proportion with them, then, however many numbers fall
between them in continued proportion, so many will also fall in
continued proportion between the numbers which have the same ratio
with the original numbers.

Proposition 9: If two numbers be prime to one another, and numbers
fall between them in continued proportion, then, however many
numbers fall between them in continued proportion, so many will also
fall between each of them and a unit in continued proportion.

Proposition 10: If numbers fall between each of two numbers and a
unit in continued proportion, however many numbers fall between each
of them and a unit in continued proportion, so many also will fall
between the numbers themselves in continued proportion.

Proposition 11: Between two square numbers there is one mean
proportional number, and the square has to the square the ratio
duplicate of that which the side has to the side.

Proposition 12: Between two cube numbers there are two mean
proportional numbers, and the cube has to the cube the ratio
triplicate of that which the side has to the side.

Proposition 13: If there be as many numbers as we please in
continued proportion, and each by multiplying itself makes some
number, the product will be proportional; and, if the original
numbers by multiplying the product make certain numbers, the latter
will also be proportional.

Proposition 14: If a square measures a square, the side will also
measure the side, and, if the side measures the side, the square
will also measure the square.

Proposition 15: If a cube number measures a cube number, the side
will also measure the side ; and, if the side measures the side, the
cube will also measure the cube.

Proposition 16: If a square number does not measure a square number,
neither will the side measure the side; and, if the side does not
measure the side, neither will the square measure the square.

Proposition 17: If a cube number does not measure a cube number,
neither will the side measure the side; and, if the side does not
measure the side, neither will the cube measure the cube.

Proposition 18: Between two similar plane numbers there is one mean
proportional number; and the plane number has to the plane number
the ratio duplicate of that which the corresponding side has to the
corresponding side.

Proposition 19: Between two similar solid numbers there fall two
mean proportional numbers; and the solid number has to the similar
solid number the ratio triplicate of that which the corresponding
side has to the corresponding side.

Proposition 20: If one mean proportional number falls between two
numbers, the numbers will be similar plane numbers.

Proposition 21: If two mean proportional numbers fall between two
numbers, the numbers will be similar solid numbers.

Proposition 22: If three numbers be in continued proportion, and the
first be square, the third will also be square.

Proposition 23: If four numbers be in continued proportion, and the
first is cube, the fourth will also be cube.

Proposition 24: If two numbers have to one another the ratio which a
square number has to a square number, and the first be square, the
second will also be square.

Proposition 25: If two numbers have to one another the ratio which a
cube number has to a cube number, and the first be cube, the second
will also be cube.

Proposition 26: Similar plane numbers have to one another the ratio
which a square number has to a square number.

Proposition 27: Similar solid numbers have to one another the ratio
which a cube number has to a cube number.

\vspace{3mm}\textbf{Book 9 Propositions }

Proposition 1: If two similar plane numbers by multiplying one
another make some number, the product will be square.

Proposition 2: If two numbers by multiplying one another make a
square number, they are similar plane numbers.

Proposition 3: If a cube number by multiplying itself, the product
will be cube.

Proposition 4: If a cube number by multiplying a cube number makes
some number, the product will be cube.

Proposition 5: If a cube number by multiplying any number makes a
cube number, the multiplied number will also be cube.

Proposition 6: If a number by multiplying itself make a cube number,
it will itself also be cube.

Proposition 7: If a composite number by multiplying any number makes
some number, the product will be solid.

Proposition 8: If as many numbers as we please beginning from a unit
be in continued proportion, the third from the unit will be square,
as will also those which successively leave out one; and the fourth
will be cube, as will also those which leave out two; and the
seventh will be at once cube and square, as will also those which
leave out five.

Proposition 9: If as many numbers as we please beginning from a unit
be in continued proportion, and the number after the unit be square,
all the rest will also be square. And, if the number after the unit
be cube, all the rest will also be cube.

Proposition 10: If as many numbers as we please beginning from a
unit be in continued proportion, and the number after the unit be
not square, neither will any other be square except the third from
the unit and all those which leave out one. And, if the number after
the unit be not cube, neither will any other be cube except the
fourth from the unit and all those which leave out two.

Proposition 11: If as many numbers as we please beginning from a
unit be in continued proportion, the less measures the greater
according to some one of the numbers which have place among the
proportional numbers.

Corollary: And it is manifest that whatever place the measuring
number has, reckoned from the unit, the same place also has the
number according to which it measures, reckoned from the number
measured, in the direction of the number before it.

Proposition 12: If as many numbers as we please beginning from a
unit be in continued proportion, by however many prime numbers the
last is measured, the next to the unit will also be measured by the
same.

Proposition 13: If as many numbers as we please beginning from a
unit be in continued proportion, and the number after the unit be
prime, the greatest will not be measured by any except those which
have a place among the proportional numbers.

Proposition 14: If a number be the least that is measured by prime
numbers, it will not be measured by any other prime number except
those the originally measuring it.

Proposition 15: If three numbers in continued proportion be the
least of those which have the same ratio with them, any two whatever
added together will be prime to the remaining number.

Proposition 16: If two numbers be prime to one another, the second
will not be to any other number as the first is to the second.

Proposition 17: If there be as many numbers as we please in
continued proportion, and the extremes of them be prime to one
another, the last will not be to any other number as the first to
the second.

Proposition 18: Given two numbers, to investigate whether it is
possible to find a third proportional to them.

Proposition 19: Given three numbers, to investigate when it is
possible to find a fourth proportional to them.

Proposition 20: Prime numbers are more than any assigned multitude
of prime numbers.

Proposition 21: If as many even numbers as we please be added
together, the whole is even.

Proposition 22: If as many odd numbers as we please be added
together, and their multitude be even, the whole will be even.

Proposition 23: If as many odd numbers as we please be added
together, and their multitude be odd, the whole will also be odd.

Proposition 24: If from an even number an even number be subtracted,
the remainder will be even.

Proposition 25: If from an even number an odd number be subtracted,
the remainder will be odd.

Proposition 26: If from an odd number an odd number be subtracted,
the remainder will be even.

Proposition 27: If from an odd number an even number be subtracted,
the remainder will be odd.

Proposition 28: If an odd number by multiplying an even number make
some number, the product will be even.

Proposition 29: If an odd number by multiplying an odd number make
some number, the product will be odd.

Proposition 30: If an odd number measures an even number, it will
also measure the half of it.

Proposition 31: If an odd number be prime to any number, it will
also be prime to the double of it.

Proposition 32: Each of the numbers which are continually doubled
beginning from a dyad is an even-times even only.

Proposition 33: If a number has its half odd, it is an even-times
odd only.

Proposition 34: If a number neither be one of those which are
continually doubled from a dyad, nor have its half odd, it is both
an even-times even and even-times odd.

Proposition 35: If as many numbers as we please be in continued
proportion, and there be subtracted from the second and the last
numbers equal to the first, then, as the excess of the second is to
the first, so will the excess of the last be to all those before it.

Proposition 36: If as many as we pleas beginning from a unit be set
out continuously in double proportion, until the sum of all becomes
prime, and if the sum multiplied into the last make some number, the
product will be perfect.

\clearpage
\end{document}